\documentclass[12pt]{article}
\usepackage{amssymb}
\usepackage{latexsym}

\newtheorem{Lemma}{Lemma}
\newtheorem{Proposition}[Lemma]{Proposition}

\newtheorem{Theorem}[Lemma]{Theorem}

\newcommand{\bck}{\!\!\!}

\newcommand{\LL}{\mbox{${\cal L}$}}
\newcommand{\CC}{\mbox{${\cal C}$}}

\def \R{\mathbb R}
\def \C{\mathbb C}
\def \L{\mathbb L}
\def \Z{\mathbb Z}
\def  \P {{\bf P}}
\def \Prob{{\bf P}}
\def \E{{\bf E}}
\def \N {{\mathbb N}}
\def \Pf {\noindent {\bf Proof.} }
\def \Re {{\rm Re}}
\def \Im {{\rm Im}}

\newcommand{\psix}[1]{(#1)^{1/2}}
\def\endpf{$\Box$}
\begin{document}
\begin{titlepage}
\begin{center}
{\large \bf The Beurling estimate for a class of random walks} \\
\end{center}
\vspace{5ex}
\begin{flushright}
Gregory F. Lawler
\footnote{Research supported in part by the
National Science Foundation DMS 9971220}$^,$\footnote{Cornell
University, Department of Mathematics, 310 Malott Hall, Ithaca, NY, 14853, lawler@math.cornell.edu}
  ~\\
Vlada Limic 
\footnote{Research supported in part by
National Science Foundation grant DMS 0104232
and NSERC Research grant}$^,$\footnote{University of British Columbia, Department of 
Mathematics, Department of Mathematics,\#121-1984 Mathematics Road,
Vancouver, BC V6T 1Z2, Canada, limic@math.ubc.ca}

\end{flushright}

\vfill

{\bf ABSTRACT:} \\ 
An estimate of  Beurling  
states that if $K$ is a curve from $0$ to
the unit circle in the complex plane,
 then the probability that a Brownian motion starting
 at $-\epsilon$ reaches the unit circle without hitting
the curve is bounded above by $ c \, \epsilon^{1/2}$.
This estimate is very useful in analysis of boundary behavior
of conformal maps,  especially for connected but
rough boundaries.  
The corresponding estimate for simple random walk was first proved
by  Kesten.     
In this note we extend this estimate 
to  random walks with zero mean, finite $(3+\delta)-$moment.
\vfill



\end{titlepage}
 \section{Introduction}

The Beurling projection theorem (see, e.g., \cite[Theorem V.4.1]{Bass})
states that if $K$ is a  closed subset of the
closed unit disk in $\C$, then the probability that a Brownian motion starting
at $-\epsilon$ avoids $K$ before reaching the unit circle is less than
or equal to the same probability for the angular projection
 \[  K' = \{|z|: z \in K \} . \]
If $K' = [0,1]$, a simple conformal mapping argument shows that
the latter probability is comparable to $\epsilon^{1/2}$ as
$\epsilon \rightarrow 0+$. 
In particular, if $K$ is a connected set of diameter one at
distance $\epsilon$ from the origin
the probability that a Brownian motion from
the origin to the unit circle avoids $K$ is bounded above by 
$c \epsilon^{1/2}.$

This estimate, which we will call the Beurling estimate,
 is very useful in analysis of boundary behavior
of conformal maps especially for connected but
rough boundaries.  
A similar estimate for random walks is useful, especially when
considering convergence of random walk to Brownian motion near
(possibly nonsmooth) boundaries.  For simple random
walk such an estimate was first  established
in \cite{Kesten} to derive a discrete harmonic measure
estimate for application to diffusion limited aggregation.
It has been used since in a number of places, e.g., in
deriving ``Makarov's Theorem'' for random walk \cite{LMakarov}
or establishing facts about intersections of random walks
(see, e.g., \cite{LPuckette}).  Recently it has been used
by the first author and collaborators to analyze the rate
of convergence of random walk to Brownian motion in 
domains with very rough boundaries.    Because of its utility,
we wish to extend this estimate to walks other than just
simple random walk.   In this note we 
extend it to a larger class of random walks.

We state the precise result in the next section, but we
will summarize briefly here.  As in \cite{Kesten},
we start with the estimate for a half-line.  We follow the
argument 
in \cite{LBook}; see \cite{M-B, Fukai}  for extensions.
The argument in \cite{LBook} strongly uses the time reversibility
of simple random walk. In fact, as was noted in \cite{Fukai}, the argument
really only needs symmetry in $x$ component. 
We give a 
 proof of this estimate,
because we need the result not just for $\Z^+$ but also for
$\kappa \Z^+$ where $\kappa$ is a positive integer.  The reason is that
we establish  the
Beurling estimate here for ``$(1/\kappa)$-dense'' sets.  One example
of such a set that is not connected
is the path of a non-nearest neighbor random walk whose increments have
finite range; a possible application of our result would be to
extend the results of \cite{LPuckette} to finite range walks.
While our argument is essentially complete for random walks that
are symmetric in the $x$ component, for the nonsymmetric case
we use a result of Fukai \cite{Fukai} that does the estimate
for $\kappa = 1$.  Since $\kappa \Z_+ \subset \Z_+$ this gives
a lower bound for our case, and our bound for the full line then
gives the upper bound.

The final section derives the general result from that
for a half-line; this argument closely follows that in
\cite{Kesten}.   
We assume a $(3 + \delta)$-moment for the increments
of the random walk in order to ensure that the asymptotics for
the potential kernel are sufficiently sharp (see (\ref{Eaexp})).  (We also use
the bound for some ``overshoot'' estimates, but in these cases weaker
bounds would suffice.)

\section{Preliminaries}

\label{Intro}
Denote by ${\Z}, {\R}, {\C}$ the integers, the real numbers and the complex numbers,
respectively.  We consider $\Z$ and $\R$ as subsets of $\C$.
Let $\Z^+ = \{k \in \Z: k > 0\}; \N = \{k \in \Z: k \geq 0\},
\Z^- =  \Z \setminus \N$.
Let $\L$ denote a discrete two-dimensional lattice (additive
subgroup) of $\C$.  Let $X_1,X_2,\ldots$ be i.i.d. random variables
taking values in $\L$  and let $S_n $ be the corresponding random
walk.  
We say that $X_1,X_2,\ldots$ generates $\L$ if for
each $z \in \L$ there is an $n$ with $\Prob(X_1 + \cdots + X_n
= z) > 0$.
Let
\[
T_B:=\inf\{l\geq 1:S_l\in B\}, \ \ T_B^0:=\inf\{l\geq 0:S_l\in B\},
\]
be the first entrance time of $B$ after time $0$, and the first entrance time of $B$ including time $0$,
respectively.
We abbreviate $T_{\{b\}},T_{\{b\}}^0$ by $T_b,T_b^0$ respectively.
Denote by $\CC_n=\{z \in \L :|z|<n\}$ the discrete open disk of radius $n$, and let
$\tau_n := T_{\CC_n^c}^0$ be the first time the random walk is not
in $\CC_n$.  

Suppose $\kappa$ is a positive integer
and $A$ is a subset of the lattice $\L$.
We
call $A$ {\em $(1/\kappa)$-dense (about the origin)}
if for every $j \in \N$, $A \cap \{j \kappa \leq
|z| < (j+1) \kappa\} \neq \emptyset$.  A set of
the form $A = \{w_j : j \in \kappa \N\}$ with
$j \leq |w_j| < j + \kappa$ for each $j$ will be called a
{\em minimal} $(1/\kappa)$-dense set.  Any
$(1/\kappa)$-dense set  contains a minimal
$(1/\kappa)$-dense set.  
If $0 < j_1 < j_2 < \infty$, 
we let $A[j_1,j_2] = A \cap (\CC_{j_2} \setminus
\CC_{j_1}).$  
If  $-\infty < j_1 < j_2 < \infty$, we
write $ [j_1,j_2]_\kappa=\kappa \N \cap [j_1,j_2)$.

The purpose of this paper is the prove the following result.
\begin{Theorem} \label{Tmain}
Suppose $\L$ is a discrete two-dimensional lattice
in $\C$ and $X_1,X_2,\ldots$ are i.i.d.
random variables that generate $\L$
such that  $\E[X_1] = 0$ and for some $\delta> 0$,
$\E[|X_1|^{3 + \delta}] < \infty$.  Then for
positive integer $\kappa$, there exists a $c < \infty$
(depending on $\kappa$ and the distribution of
$X_1$) such
that for every $(1/\kappa)$-dense set $A$ and every
$0 < k < n < \infty$,
\[ \P(\tau_{2n} < T_{A[k,n]})   \leq c \sqrt{k/n} .
\]
\end{Theorem}

We start by making some reductions.  
Since $B\subset A$ clearly implies
$\Prob(\tau_m < T_A^0)\leq \P(\tau_m < T_B^0)$, it suffices to prove the theorem
for minimal $(1/\kappa)$-dense
sets  $A = \{w_j: j \in \kappa \N\}$ and, without loss of generality,
we assume that $A$ is of this form.
By taking a linear
transformation of the lattice if necessary, we
may assume that 
 ${\L}$ is of the form
\[ \L=\{j+kz^*:j,k\in \Z\},
\]
where $z^* \in \C \setminus \R$ and that the
covariance matrix of $X_1$ is a multiple of
the identity.  (When dealing with mean zero, finite
variance lattice random walks, one can always choose the
lattice to be the integer lattice in which case one may have
a non-diagonal covariance matrix, or one can choose a more
general lattice but require the covariance matrix to be
a multiple of the identity.  We are choosing the latter.)
Let $p$ be the (discrete)
probability mass function of $X_1.$  Then our assumptions
are
 $\{z:p(z) > 0 \}$ generates $\L$ and   for some $\delta,\sigma^2 > 0$, 
\begin{equation}  \label{cond1}
  \sum_{z} z p(z) = 0,
\end{equation} 
\begin{equation}  \label{cond2}
 \sum_{z}  \Re(z)^2 p(z)
  = \sum_z \Im(z)^2 p(z)  = \sigma^2 > 0 , 
\end{equation}
\begin{equation}  \label{cond3}
  \sum_{z} |z|^{3+\delta} p(z) < \infty ,
\end{equation}
Let $p_*(z) =p(z)$ be step probability mass function of the time-reversed walk; 
and note that $p_*$ also satisfies (\ref{cond1})-(\ref{cond3}).
We denote by $\P_*(A)$ the probability of $A$ under steps according to $p_*$. 
 We call a function $f$
{\em $p$-harmonic} at $w$ if
\begin{equation}
\label{Eharm}
              \Delta_p f(w) := \sum_{z} p(z) \, [f(z+w) - f(w)] = 0 . 
\end{equation}
Let $X_1,X_2,\ldots$ be independent $\L$-valued
random variables with
  probability mass function $p$,
and let  $S_n=S_0+\sum_{i=1}^n X_i, n\geq 0$ be the corresponding
random walk.   
Denote by $\P^x$ (resp., $\E^x$) the law (resp., expectation)
of $(S_n,n\geq 0)$ when $S_0=x$, and we will write $\Prob,\E$,
for $\P^0,\E^0$.   

\medskip

Let $a(z)$ denote the {\em potential kernel} for $p$,
\[
a(z) =\lim_{n\to \infty} \sum_{j=0}^n \left[\Prob(S_j=0)-\Prob(S_j=z)\right],
\]
and let $a^*(z)$ denote the potential kernel using $p_*$.
Note that $a$ is $p_*$-harmonic and $a^*$ is $p$-harmonic for $z \neq 0$
and $\Delta_{p_*} a(0) = \Delta_pa^*(0) = 1$.
In \cite{FU} it is shown that 
under the assumptions (\ref{cond1}) - (\ref{cond3})
there exist constants $\bar k,c$ (these constants, like all constants
in this paper, may depend on $p$), such that for all
$z$,
\begin{equation}
\label{Eaexp}
\left|\, a(z) - \frac{\log |z|}{\pi \, \sigma^2 } - \bar{k}\, \right|  \leq \frac{c}{|z|}.
\end{equation}\
Since $a^*(z) = a(-z)$, this also holds for $a^*$.

%
%

As mentioned above,  $\CC_n=\{z \in \L :|z|<n\}$ is the discrete open disk of radius $n$ and
$\tau_n := T_{\CC_n^c}^0$. 
Denote by $\LL_n$ the discrete open strip $\{x+iy \in \L: |y|< n\}$ of width $2n$ and
let
$ \rho_n:= T_{\LL_n^c}^0,$  i.e., $\tau_n,\rho_n$ are 
the exit times from the disk and the strip, respectively.

For any proper subset $B$  of $\L$, let $G_B(w,z)$
denote the  Green's function of  $B$
defined by
\begin{equation}
\label{EgreenB}
G_B(w,z) = \sum_{j=0}^{T^0_{B^c}-1} \Prob^w(S_j = z).
\end{equation}
This equals zero unless $w,z\in B$.  
We will write $G_n$ for $G_{\CC_n}$.
If $w,z\in B$, and 
$G(w) := G_B(w,z)$, then $\Delta_pG(w) = - \delta(w-z)$ 
where $\Delta_p$ is as in (\ref{Eharm}), and where $\delta(\cdot)$ is the Kronecker symbol
$\delta(x)=1,\, x=0$ and $\delta(x)=0,\,x\neq 0$.
  Let $G_B^*(w,z) $ denote the Green's function
for $p_*$ and note that $G_B(z,w) = G_B^*(w,z)$.  
A useful formula for finite $B$ is
\begin{equation}  \label{Euseful}
 G_B(w,z) = \E^w[a^*(S_{T} - z)] - a^*(w-z) 
   = \E^w[a(z - S_T)] - a(z-w), 
\end{equation}
where $T = T_{B^c}^0$.  This is easily verified by noting that for
fixed $z \in B$, each of the three expressions describes the function
$(w)$ satisfying: $f(w) = 0, w \not\in B; \Delta_p f(z) = -1;
\Delta_pf(w) = 0, w \in B \setminus \{z\}$.  The following
``last-exit decomposition''
 relates the Green's function and escape probabilities:
\begin{equation}  \label{lastexit}
 \Prob^z\{T^0_{B'} < T^0_{B^c}\} =
              \sum_{w \in B'} G_B(z,w) \, \Prob^w(T_{B'} > T_{B^c}).
\end{equation}
It is easily derived by focusing on the last visit to $B'$ strictly less
than $T^0_{B^c}$.  
 

  
For the remainder of this paper we fix $p,\kappa$ and
allow constants to depend on $p,\kappa$.   We assume $k \leq n/2$,
for otherwise the inequality is immediate.
The values of universal constant may change from line to line without further notice.
In the next two sections will prove that
\[ \P(\tau_{n^4} < T_{[k,n]}) \leq \frac{1}{\log{n}} \sqrt{\frac{k}{n}}. \]
(Here, and throughout this paper, we use $\asymp$ to mean that both
sides are bounded by constants times the other side where the constants
may depend on $p,\kappa.$)
In the final section we establish the uniform upper bound 
for all minimal $(1/\kappa)$-dense sets.

\section{Green's function estimates}

We start with an ``overshoot'' estimate.  

\begin{Lemma} \label{overshoot} There is a $c$ such that for all $n$ and all $z$ with $|z| < n$,
\[   \E^z[|S_{\tau_n}| ]  \leq n + c\, {n^{2/3}} , \;\;\;\;
 \E^z[\log |S_{\tau_n}| - \log n] \leq c \, n^{-1/3} . \]
\end{Lemma}

\Pf If $a > 0$, 
since $\{|S_{\tau_n}|-n\geq a\}\subset\{|X_{\tau_n}|\geq a\}$ we have
\[
\Prob^z(|S_{\tau_n}|-n\geq a) \leq
\sum_{j=1}^\infty \Prob^z(\tau_n =j,|X_j|\geq a)\leq
 \sum_{j=1}^\infty \Prob^z(\tau_n >j-1,|X_j|\geq a)\]
\[  \hspace{1in}  \leq \sum_{j=1}^\infty \Prob^z(\tau_n > j-1) \, \Prob(|X_1|
    \geq a )    \leq \E^z(\tau_n) \, \Prob(|X_1| \geq a ) . \]  From the
central limit theorem, we know that $\Prob^z\{\tau_n > r + n^2 \mid \tau_n > r \}
< \alpha < 1$.  Therefore, 
$\tau_n/n^2$ is stochastically bounded by a geometric random variable
with success probability $1-\alpha$,
and hence $\E^z[\tau_n] \leq c n^2 . $  Since $\E[|X_1|^{3}]
< \infty$,
\begin{equation}  \label{overshoot1}
  \Prob(|X_1| \geq b ) = \Prob(|X_1|^{3} \geq
    b^{3}) \leq c \, b^{-3 } . 
\end{equation}
Therefore,
\begin{equation}  \label{overshoot2}
   \Prob^z(|S_{\tau_n}|-n\geq a\, n^{2/3} ) \leq c \,  a^{-3}, 
\end{equation}
and \[ 
\     \E[|S_{\tau_n}| - n] \leq n^{2/3} + \int_{n^{2/3}}
  ^\infty \Prob^z (|S_{\tau_n}|-n\geq y ) \; dy \leq c \, n^{2/3} . \]
The second inequality
follows immediately using $\log(1+x) \leq x$.

\medskip


\noindent {\bf Remark.}  With a finer argument, we could show, in fact,
that $\E^z[|S_{\tau_n}|] \leq n +c$.
  By doing the more refined estimate we
could improve some of the propositions below, e.g., the $O(n^{-1/3})$
error
term in the next proposition is actually $O(n^{-1})$.  However, since the
error terms we have proved here suffices for this paper, we will not prove
the sharper estimates.

\begin{Lemma}\label{Greenlemma1}
\[  \pi \, \sigma^2 \, G_n(0,0) = \log n +O(1) . \]
If $|z| < n$,
\[   \pi \, \sigma^2 G_n(0,z) = \log n - \log |z| + 
   O\left(\frac{1}{|z|}\right) +O(n^{-1/3}), \]
\[   \pi \, \sigma^2 G_n(z,0) = \log n - \log |z| + 
   O\left(\frac{1}{|z|}\right) +O(n^{-1/3}). \]
Also, for every $b < 1$, there exist $c > 0$ and
$N$ such that for all $n \geq N$,
\begin{equation}  \label{jun22.1}
       G_n(z,w) \geq c , \;\;\;\;  z,w \in \CC_{bn} . 
\end{equation}
\end{Lemma}

\Pf The first expression  follows from (\ref {Eaexp}), (\ref{Euseful}) and
Lemma \ref{overshoot} since
$a(0) = 0$.  The next two expressions again use
(\ref{Euseful}), Lemma \ref{overshoot}, and (\ref {Eaexp}).
For the final expression, first note it is true for $b = 1/4$, since
for $0 \leq |z|,|w| < n/4$,
 $G_n(z,w) \geq G_{3n/4}(0,w-z)$. 
For $b < 1$, the invariance principle implies that there
is a $q = q_b > 0$ such that for all $n$ sufficiently large,
with probability at least $q$ the random walk (and reversed
random walk) starting at $|z| < bn$ reaches $\CC_{n/4}$
before leaving $\CC_n$.   Hence, by the strong Markov
property, if $|z| < bn, |w| < bn$,
$G_n(z,w) \geq q \, \inf_{|z'| < n/4} G_n(z',w).$  Similarly,
using the reversed random walk, if $|w| < bn, |z'| < n/4$,
$G_n(z',w) \geq q \, \inf_{|w'| < n/4} G_n(z',w')$.
\ \ \ \endpf

\begin{Lemma}\label{Greenlemma}
  If  $m \geq n^4$ and $|z|,|w| \leq n$, 
\[  \pi \, \sigma^2 \, G_m(z,w) =
\log m - \log|z-w| + 
 O\left(\frac{1}{|z-w|}\right) +O(n^{-4/3}).
    \]
\end{Lemma}

\Pf    Since $G_{m-n}(0,w-z) \leq G_m(z,w) \leq
G_{m+n}(0,w-z)$, this follows from the previous
lemma.

\begin{Lemma}
There is a $c<\infty$ such that for every $z \in \CC_n$
and every minimal $(1/\kappa)$-dense set $A$,
\begin{equation}  \label{jun22.2}
  \sum_{w \in A} G_n(z,w) \leq c \, n . 
\end{equation}
\end{Lemma}

\Pf By Lemma \ref{Greenlemma1},
\[ \pi \, \sigma^2 \,
G_n(z,w)\leq \pi \, \sigma^2 \, G_{2n}(0,w-z) \leq 
     \log n - \log |w-z| + O(1) . \]
If $A$ is a minimal $(1/\kappa)$-dense set, then
$  \# \{w \in A: |z-w| \leq r\}  \leq cr$, for some $c$ independent
of $z$.  Hence,
\[ \sum_{w \in A} G_n(z,w) \leq c \, \sum_{j=1}^{2n}
    [\log n - \log j + O(1) ] = O(n) . \] \endpf

\section{Escape probability estimates for $[j,k]_\kappa$}
\label{EPE}


The main purpose of this section is to obtain estimates in Proposition \ref{prop2} and Lemmas \ref{loglem}
and \ref{LPesc} which will be used in the proof of Theorem \ref{Tmain} in section \ref{PT}.
 
\begin{Lemma} 
\begin{equation}  \label{jun20.2}
   \P(\tau_n< T_{[-n,n]_\kappa})\asymp
\frac 1n . 
\end{equation}
\end{Lemma}

\Pf  Let $q(n) =  \P(\tau_n\leq T_{[-n,n]_\kappa})$ and note
that if $k \in [-n/2,n/2]_\kappa$, then
\[        q(4n) \leq              \P^k(\tau_n\leq T_{[-n,n]_\kappa})
                              \leq q(n/4) . \]
The last-exit decomposition (\ref{lastexit}) tells us
\[                   \sum_{k \in [-n/2,n/2]_\kappa}
                    G_n(0,k) \,  \P^k(\tau_n< T_{[-n,n]_\kappa}) \leq \sum_{k \in [-n,n]_\kappa}
                    G_n(0,k) \,  \P^k(\tau_n< T_{[-n,n]_\kappa}) =1 . \]
But (\ref{jun22.1}) and (\ref{jun22.2}) imply that
\[           \sum_{k \in [-n/2,n/2]_\kappa}
                     G_n(0,k)  \asymp n , \ \mbox{ \endpf}\]
which gives $q(4n)=O(1/n)$.
The lower bound can be obtained by noting
$\P(\rho_n < T_{\Z})\leq \P(\tau_n< T_{[-n,n]_\kappa})$  
which reduces the estimate to a  one-dimensional
``gambler's ruin'' estimate in the $y$-component.  This can be established in a number
of ways, e.g., using a martingale argument.
\ \ \ \endpf

\begin{Lemma}  \label{gottahit}
There exist $c >0$ and $N < \infty$ such that
if $n \geq N$ and $z \in \CC_{3n/4}$,
\[  \Prob^z(T^0_{A[n/4,n]} < \tau_n)
 \geq\Prob^z(T^0_{A[n/4,n/2]} < \tau_n) \geq c . \]
\end{Lemma}

\Pf  Let 
\[                  V = \sum_{j=0}^{\tau_n - 1} 1\{S_j  \in A[n/4,n/2]\} , \]
be the number of visits to $A[n/4,n/2]$ before leaving
$\CC_n$.  Then (\ref{jun22.1}) and (\ref{jun22.2}) show
that there exist $c_1,c_2$ such that for $n$ sufficiently large,
\[   c_1 \, n \leq \E^w[V] \leq c_2 \, n , \;\;\;\;  w \in \CC_{7n/8} . \]
In particular, if $z \in \CC_{3n/4}$, 
\[ c_1 \, n \leq \E^z[V] = \Prob^z(V \geq 1)\;  \E^z[V \mid V \geq 1]
                  \leq c_2 \, n \, \Prob^z(V \geq 1) .\ \mbox{ \endpf} \]

\begin{Lemma}  \label{logest}
There exist  $0 < c_1 < c_2 < \infty$ and $N < \infty$
such that if $n \geq N$,  
\[       \frac {c_1}{\log n} \leq \Prob^z(T_0 < \tau_n) \leq \frac{c_2}{\log n},
\;\;\;\;  z \in \CC_{9n/10} \setminus \CC_{n/10} . \]
\end{Lemma}

\Pf   This
follows immediately from Lemma \ref{Greenlemma1} and 
$ G_n(z,0) = \Prob^z(T_0 < \tau_n) \, G_n(0,0)$. \ \endpf\\

Let $T^+ = T_{\kappa \Z^+}, T^- = T_{\kappa \Z \setminus \Z^+}  .$
Define  
\[
E_n^+ =\{\rho_n < T^+\},\
E_n^- =\{\rho_n < T^-\},\ \tilde{E}_n^- =\{\rho_n < T_{\kappa \Z^-}\}
\]
and
\[
E_n =E_n^+\cap E_n^-=\{\rho_n <T_{\kappa \Z}\} .
\]
Recall that $\P_*$ stands for the probability under step distribution $p_*$.
\begin{Lemma}
\label{Lindep}
$\P(E_n) = \P_*(\tilde{E}_n^-) \P(E_n^-)$.
\end{Lemma}

\Pf 
Consider $E_n^-\cap (E_n^+)^c =V_1\cup V_2 \cup \ldots $,
where
\[
V_m=\{\rho_n < T_{\{\ldots,-\kappa,0,\kappa,\ldots,\kappa( m-1)\}}\} \cap
\{\rho_n > T_{\kappa m}\},
\]
is the event that integer $\kappa m$ is the smallest integer
in $\kappa \Z$  visited
by the walk before time $\rho_n$.
Clearly $V_1,V_2,\ldots$ are disjoint events.  Write
\[
V_m =\bigcup_{j=1}^\infty V_{m,j},
\]
where $V_{m,j} :=V_m \cap \{S_j=\kappa m\} \cap \{S_l \neq \kappa m, l=j+1,\ldots,\rho_n\}$
is the intersection of $V_m$
with the event that  $\kappa m$ is visited for the last time (before time $\rho_n$)
at time $j$.
Again, $V_{m,j}$ are mutually disjoint events.
Therefore,
\begin{equation}
\label{Esta}
\Prob(E_n^-\cap (E_n^+)^c) =\sum_{m=1}^\infty \sum_{j=1} ^\infty \Prob(V_{m,j}).
\end{equation}
Note that due to the strong Markov property, and homogeneity of the line
and the lattice, we have
\begin{eqnarray}
\bck \Prob(V_{m,j}) \bck&=&\bck \Prob(S_j=\kappa m, j< \rho_n\wedge 
T_{\kappa\{\ldots,-1,0,1,\ldots,m-1\}}) \, \P^{\kappa m}
  (\rho_n < T_{\kappa \{\ldots,-1,0,1,\ldots,m\}})\nonumber\\
\label{EVmj} \bck&=&\bck \Prob (S_j=\kappa m, j-1<\rho_n\wedge
 T_{\kappa\{\ldots,-1,0,1,\ldots,m-1\}} ) \, \P(E_n^-).
\end{eqnarray}
By reversing the path we can see that 
\[ \Prob(S_j=\kappa m, j-1<\rho_n\wedge T_{\kappa\{\ldots,-1,0,1,\ldots,m-1\}} )\hspace
{2in}  \]
\begin{equation} \hspace{2in}
\label{Esym} =
\P_*^{\kappa m}(S_j=0, j-1<\rho_n\wedge T_{\kappa\{\ldots,-1,0,1,\ldots,m-1\}} ).
\end{equation}
Also note that
\[ \P_*^{\kappa m}(S_j=0, j-1<\rho_n\wedge T_{\kappa\{\ldots,-1,0,1,\ldots,m-1\} })
=  \hspace{2in}  \]
\begin{equation}
\label{Etra}  \hspace{2in}
\P_*(S_j=-\kappa m, j-1<\rho_n\wedge T_{\kappa\{\ldots,-2,-1\}} )
\end{equation}
by translation invariance.
Now,
\[
\{S_j=-\kappa m, j-1<\rho_n\wedge T_{\kappa\{\ldots,-2,-1\}}\}=
\{\rho_n\wedge T_{\kappa\Z^-} = T_{-\kappa m} = j\}
\]
and since
\[
\sum_{m=1}^\infty \sum_{j =1}^\infty \P_*(\rho_n\wedge T_{\kappa\Z^-} = T_{-\kappa m} = j)=
\sum_{m=1}^\infty \P_*(\rho_n\wedge T_{\kappa\Z^-} = T_{-\kappa m})=\P_*(\rho_n> T_{\kappa\Z^-})
\]
relations (\ref{EVmj})-(\ref{Etra}) imply
\[
\sum_{m\geq 1}\sum_{j\geq 1} \P(V_{m,j}) =\P_*(\rho_n> T_{\kappa\Z^-})\,
\P(E_n^-)=\P_*((\tilde{E}_n^-)^c)\, \P(E_n^-).
\]
This together with (\ref{Esta}) implies the lemma.
\ \ \ \endpf

\medskip

\noindent{\bf Remark.}
The above result implies the following remarkable claim: if the step distribution of the walk 
is symmetric with respect to $y$-axis then, under $\P$, the events
$E_n^+$ and $E_n^-$ are independent.
 
\medskip

\noindent {\bf Remark.}  Versions of this lemma have appeared in a number of
places. 
See \cite{LBook,M-B,Fukai}.


\begin{Lemma}
\label{Cabc}
\begin{equation}  \label{Esecond}
\P(\rho_n\leq T_{\kappa\N} )\asymp \P(\rho_n\leq T_{\kappa\Z^-})
  \asymp \frac{1}{\sqrt{n}}.
\end{equation}
\end{Lemma}
\Pf In the case $\kappa =1$, this was essentially proved by Fukai
\cite{Fukai}.  Theorem 1.1 in \cite{Fukai} states that
\begin{equation}
\label{ETfukai}
\Prob(n^2 < T_{\N})\asymp \frac{1}{n^{1/2}}.
\end{equation}
for any zero-mean
aperiodic random walk on lattice $\Z^2$ with $2+\delta$ 
finite moment. 
Note that we can linearly map $\L$ onto $\Z^2$, and by this cause only 
multiplicative constant change (depending on $\L$) 
in the conditions (\ref{cond1})-(\ref{cond3}), which imply the assumptions 
needed for (\ref{ETfukai}) to hold.  The conversion from $n^2$ to $\rho_n$
is not difficult and his argument can be extended to give this.  Note that
this gives a lower bound for other $\kappa$,  
%
\begin{equation}
\label{Efukai}
\P(\tau_n < T_{\kappa\N})\geq \frac{c}{n^{1/2}},
\end{equation}
where $c$ depends on $\L$ and transition probability $p$ only.
Hence, the two terms in the product in Lemma \ref{Lindep} are bounded below by $c/\sqrt n$
but the product is bounded above by $c_1/n$.  Hence, each of the terms is also bounded
above by $\tilde c/\sqrt n$, 
 and this proves the statement.
\ \ \ \endpf

\medskip

\begin{Lemma}
\label{Lhelp}
There exists $c\in (0,\infty)$ such that\\
(a) $P^{-n}(T_{-n}<T_{\kappa \N})\leq 1-\frac{c}{\log{n}} $,\\
(b) If $|z|\geq n $ then $P^{z}(T_{z}<T_{\kappa \Z})\leq 1-\frac{c}{\log{n}} $.
\end{Lemma}

\Pf 
We prove (a), and note that (b) can be done similarly.
It is equivalent to show
\[
\Prob(T_{\kappa\N+n}<T_{0})\geq \frac{c}{\log{n}}
\]
Note that since $\tau_n \leq T_{\kappa \N+n}$, Lemma \ref{Greenlemma1} 
yields the upper bound on the above probability 
of the same order.
For the lower bound note that invariance principle implies
\begin{equation}
\label{Einva}
\Prob(\tau_n<T_0,\Re(S_{\tau_n})\geq 4n/5)\geq  \frac{\Prob(\tau_n<T_0)}{100}\geq \frac{c}{\log{n}},
\end{equation}
by Lemma \ref{Greenlemma1}.
Use Markov property and Lemma \ref{gottahit} applied to disk centered at $n=(n,0)$ of radius $9n/10$ to get 
\[
\Prob(T_{\kappa \N+n}<T_0|\tau_n<T_0,\Re(S_{\tau_n})\geq 4n/5,|S_{\tau_n}|-n\leq n/5)\geq c,
\]
uniformly in $n$.
An easy overshoot argument yields
$\P(\tau_n<T_0,\Re(S_{\tau_n})\geq 4n/5, |S_{\tau_n}|-n\leq n/5)\asymp \P(\tau_n<T_0,\Re(S_{\tau_n})\geq 4n/5)$,
which implies the lemma. 
\ \ \ \endpf

\medskip
\begin{Proposition}  \label{prop2}
If $j,n\in \Z^+$,
\label{Pabcd}
\[ (a)\;\;  \P(\tau_n\leq T_{\kappa \N}) \asymp \P(\tau_n \leq
   T_{\kappa \Z^+}) 
 \asymp \frac{1}{\sqrt{ n}}, \]
\[  (b) \;\;
 \Prob^{-n}(S_{T_{\kappa \N}} = 0 ) =O\left(\frac{1}{\sqrt{ n}}\right), \]
\[ (c)\;\; \Prob^{n}(S_{T_{\kappa \N}} = 0) =O\left(\frac 1 {n^{3/2}}\right).\]
\[ (d) \;\;\P(\tau_n< T_{\kappa(j+\N)})=O\left(\sqrt{\frac{j}{n}}\right),\]
\[ (e) \;\; \P(\tau_n< T_{ \kappa(-j +  \N)})=O\left( \frac{1}{\sqrt{ jn}}\right),\]
\end{Proposition}

\Pf $(a) \; $  
A simple Markov argument gives
\[  \Prob(\tau_n \leq T_{\kappa \N}) \leq
 \Prob(\tau_n \leq T_{\kappa \Z_+}) \leq
\Prob(S_{T_{\kappa \N} }\neq 0)^{-1} \;
  \Prob(\tau_n \leq T_{\kappa \N}), \]
and hence the first two quantities are comparable.
Since $\tau_n \leq \rho_n$, 
(\ref{Esecond})  gives  $\P(\tau_n\leq T_{\kappa \N}) \geq c/
\sqrt n$.  For the upper bound,  let $A^- = A^-_n$ be the event
that $\Re(S_{\tau_n}) \leq 0$.  By invariance principle, $\Prob(A^-) \geq
1/4$.  However, we claim that $\Prob(A^- \mid \tau_n\leq T_{\kappa \N})
\geq \Prob(A^-)$.  Indeed, by translation invariance, we can see for
every $j > 0$, $\Prob^{j\kappa}(A^-) \leq \Prob(A^-)$, and hence by
the Strong Markov property, $\Prob(A^- \mid \tau_n >  T_{\kappa \N})
\leq \Prob(A^-)$.  Therefore,
\[  \P(\tau_n \leq T_{\kappa \N}, \Re(S_{\tau_n}) \leq 0) \geq (1/4)
          \,  \P(\tau_n \leq T_{\kappa \N}). \]
The invariance principle can now be used to see that for some $c$,
\[  \P(\rho_n \leq  T_{\kappa \N} \mid 
\tau_n \leq T_{\kappa \N}, \Re(S_{\tau_n}) \leq 0) \geq c , \]
and hence $\P(\rho_n \leq  T_{\kappa \N}) \geq (c/4) \,  \P(\tau_n \leq T_{\kappa \N}).$

$(b) \;$   Let $T = T_{-n} \wedge T_{\kappa \N}$.  
Since $\Prob^{-n}(S_T \neq -n) \geq c/\log n$ by Lemma \ref{Lhelp}(a), it suffices by the strong
Markov property to show that 
\[               \Prob^{-n} (S_T = 0) \leq  \frac{c}{(\log n) \, \sqrt n}. \]
By considering reversed paths,  
we see that
\[             \Prob^{-n}(S_T = 0) =\Prob_*(S_T =-n) . \]
But 
\begin{eqnarray*}
\Prob_*(S_T = -n) & = & \Prob_*(\tau_{n/2} < T_{\kappa \N}) \, \Prob_*(S_T = -n
\mid \tau_{n/2} < T_{\kappa \N}) \\
 & \leq & \Prob_*(\tau_{n/2} < T_{\kappa \N}) \,
            \Prob_*(S_T = -n
\mid \tau_{n/2} < T_{\kappa \N}, |S(\tau_{n/2})| \leq 3n/4)\\
  && \hspace{2in}+ \Prob_*(|S(\tau_{n/2})| \geq 3n/4) .
\end{eqnarray*}
To bound the last line, note that
by (\ref{Esecond}), $\Prob_*(\tau_{n/2} < T_{\kappa \N})
\leq c /\sqrt n$ and the conditional probability is bounded
by a term of order $1/\log n$ due to Lemmas \ref{gottahit} and 
\ref{logest}.
Inequality (\ref{overshoot2}) implies that
$\Prob_*(|S(\tau_{n/2})| \geq 3n/4) \leq c/n$.

$(c)$  We will start with the estimate
\begin{equation}  \label{jun23.1}
 \Prob^z(S_{T_{\kappa \Z}} = w) \leq \frac{c}{n} 
\mbox{ if } |z-w| \geq n. 
\end{equation}
Without loss of generality assume $w=0$, $|z| \geq n$.
As in (b), it suffices to show that  $\Prob^z(S_{ T_z \wedge T_{\kappa \Z}} = 0)
\leq c/(n\log n)$ due to Lemma \ref{Lhelp}(b).
By using reversed paths, 
we see that $\Prob^z(S_{ T_z \wedge T_{\kappa \Z}} = 0)  = \Prob_*
(S_{T_{\bar z} \wedge T_{\kappa \Z}} = \bar z) .$
Hence
it suffices to show that for all $|z| \geq n$,
\[         \Prob_*(S_T = z ) \leq \frac{c} {n \log n} , \]
where $T = T_z \wedge T_{\kappa \Z}$.   Similarly to
(b), we have $\Prob_*(\tau_{n/2} < T_{\kappa \Z}) \leq n^{-1}$ and
$\Prob_*(S_T = z \mid \tau_{n/2} < T_{\kappa \Z},
   |S_{\tau_{n/2}}| \leq 3n/4) \leq c / \log n$ .  
We have to
be a little more careful with the second term, but
\begin{eqnarray}
\lefteqn{ \Prob_*(|S_{\tau_{n/2}}| \geq 3n/4, \tau_{n/2} < T_{\kappa \Z}) }\nonumber\\
& \leq & \P_*(|S_{\tau_{\sqrt n}}| \geq n/2) +
                \P_*(\tau_{\sqrt n} < \tau_{n/2} \wedge T_{\kappa \Z},
               |S_{\tau_{n/2}} | \geq 3n/4\}\nonumber\\
& \leq &  O(n^{-2}) + O(n^{-1/2})O(n^{-1}) = O(n^{-3/2}). \label{Ecareful}
\end{eqnarray}
Using (b) and (\ref{jun23.1}) and noting $\{S_{T_{\kappa \N}}=0\}=\cap_{k=0}^\infty
\{S_{T_{\kappa \Z}}=-k\}\cap\{S_{T_{\kappa \N}}\circ \theta_{T_{\kappa \Z}}=0\}$, we conclude that
\begin{equation}
\label{Econclude}
\Prob^z(S_{\kappa \N} = 0) \leq \frac{c}{|z|^{1/2}} , \;\;\;\;
  |z| \geq n. 
\end{equation}
The remainder of the argument
 is done similarly to (b).   
Namely,
use estimate (\ref{Ecareful}) and note that
the probability that the random
walk starting at $n$ reaches a distance of $n/2$ from its starting
point without hitting $\kappa \N$ is $O(n^{-1})$, and, given that $|S_{\tau_{n/2}}|\leq 3n/4$,
the probability that it afterwards enters $\kappa \N$ at the origin
is $O(n^{-1/2})$ due to (\ref{Econclude}).

$(d) \;$  We may assume $j \kappa \leq n/4$.   By the Markov
property, translation invariance, (a), and (b), 
if $l \kappa \leq n/4$,
\begin{eqnarray*}
\lefteqn{\P(\tau_n< T_{\kappa(l+1+ \N)}) - \P(\tau_n< T_{\kappa(l+ \N)})}
\hspace{1.6in} \\
& = & \Prob(T_{l \kappa} < T_{\kappa(l+1 + \N)} \wedge \tau_n)
  \, \Prob^{l\kappa}(\tau_n < T_{\kappa(l+1+ \N) }) \\
 &\leq & \P(T_{l \kappa}=  T_{\kappa(l + \N)}) \, \Prob(\tau_{n/2} < T_{\kappa \Z^+})\\
 &= & \P^{-l \kappa}(S_{\kappa \N} = 0) \, \Prob(\tau_{n/2} < T_{\kappa \Z^+})\\
    & \leq &  c/\sqrt {ln}
\end{eqnarray*}
If we sum this estimate over $l=0,\ldots,j$, we get the estimate.

$(e)\;$ This is done similarly to (d), using (c) instead of (b).   
\ \ \ \endpf

\begin{Lemma}  \label{loglem}
There exist $0 < c_1 < c_2 < \infty$
and $N < \infty$, such that if $n \geq N, m=n^4$, and\\
(i) if $w \in \CC_{4n} \setminus \CC_{3n}$,
\[  \frac{c_1}{\log n} 
                \leq \Prob^w(\tau_m < \eta_n ) \leq  \frac{c_2}{\log n} , \]
where $ \eta_n = \inf\{j: |S_j| \leq 2n\}$.\\
(ii) if $w \in \CC_{4n}$,
\[  \Prob^w(\tau_m < T_{A[n/2,n]}) \leq  \frac{c_2}{\log n} . \]
\end{Lemma}

\noindent
{\bf Remark.} When $w \in \CC_{4n} \setminus \CC_{3n}$   (i) implies
a lower bound of the same order in (ii).

\medskip

\Pf  
(i) Let $T = \tau_m \wedge \eta_n$.   We will show that
\begin{equation}  \label{jun23.2}
\P^w(T = \tau_m) \asymp 1/\log n.  
\end{equation}
  Consider the martingale
$M_j = \pi \, \sigma^2 \,
[a^*(S_{j \wedge T}) - \bar k] - \log n $, and note
that $M_j = \log |S_{j \wedge T}| - \log n +O(|S_{j \wedge T}|^{-1})$.
Therefore,
\begin{equation}  \label{jun23.3}
    \log 3 + O(n^{-1}) \leq M_0 \leq \log 4 + O(n^{-1}) . 
\end{equation}
The optional sampling theorem implies that 
\[  \E^w[M_0] = \E^w[M_T] = \hspace{3in} \]
\begin{equation}  \label{jun23.4}
 \hspace{.6in} \E^w[M_T 1_{\{|S_T| \geq m\}}] + \E^w[M_T1_{\{|S_T|<n\}}] + 
\E^w[M_T1_{\{ n \leq |S_T| \leq  2n \}} ] 
\end{equation}
(the
estimate (\ref{overshoot2}) can be used to show that the optional
sampling theorem is valid).
Note that
\begin{eqnarray*} 
(\log m)  \, \Prob^w(T = \tau_m) & \leq& 
\E^w[\log |S_T|1_{\{ T  = \tau_m \}}  ]  \\  & \leq  &
(\log m ) \, \Prob^w(T = \tau_m) + \E^w[\log |S_{\tau_m}| - \log m]\\
&          \leq  & 4\, (\log n) \, \Prob^w(T = \tau_m) + O(n^{-4/3}) . 
\end{eqnarray*}
The last inequality uses Lemma \ref{overshoot}.  Therefore,
\[  \E^w[M_T1_{\{ |S_T| \geq m \}}] = 3 \, (\log n) \, \Prob^w(T = \tau_m)
                  +O(n^{-4/3}), \]
and hence it suffices to show that 
\begin{equation}
\label{EMorder1}
\E^w[M_T1_{\{|S_T| \geq m\}}]  \asymp 1.
\end{equation}
Clearly,
\[   O(n^{-1}) \leq \E^w[M_T1_{\{  n \leq |S_T| \leq 2n \}}] \leq
            \log 2 + O(n^{-1}) . \]
Also,
\begin{eqnarray*}
\P^w(|S_T| < n) & = & \sum_{|z| < n} \sum_{w' \in \CC_m \setminus \CC_{2n}}
                G_{\CC_m \setminus \CC_{2n}}(w,w') \, \P(X_1 = z - w') \\
     & \leq & c \, (\log n) \,  \sum_{w' \in \CC_m \setminus \CC_{2n}} \sum_{|z| < n}
                    \P(X_1 = w'-z) \\
  & \leq & c \, n^2 \, (\log n) \sum_{|z'| \geq n} \Prob(X_1 = z_1') \\
   & \leq & c \, n^{-1}\,  (\log n)\, \E[|X_1|^3]\,  \leq \, c \, n^{-1}\,  \log n, \\
       \end{eqnarray*}
and hence $\E^w[M_T; |S_T| < n] = O(\log^2 n/n)$.  
Combining these estimates with (\ref{jun23.3}) and (\ref{jun23.4})
gives (\ref{EMorder1}) and therefore (\ref{jun23.2}).

(ii) Let $q = q(n,A)$ be the maximum of $ \Prob^w(\tau_m < T_{A[n/2,n]})$
where the maximum is over all $w \in \CC_{4n}$.  Let $w = w_n$ be a point
obtaining this maximum.  Let $\bar \eta_n$ be the first time that a random walk
enters $\CC_n$ and let $\eta^* _n$ be the first time after this time that 
the walk leaves $\CC_{2n}$.  Then by a Markovian argument and an
easy overshoot argument we get
\[  \P^z(\tau_m < T_{A[n/2,n]}; \bar \eta_n < \tau_m) \leq
          \alpha \, q +O(n^{-1}) , z\in \CC_{4n}\]
where $\alpha=1-c < 1$ for $c$ the constant from Lemma \ref{gottahit}.    The
$O(n^{-1})$ error term comes from considering the probability that
$|S_{\eta^*_n}|  \geq 4n$.   
By letting $z=w$ we get
\[ \P^w(\tau_m < T_{A[n/2,n]}) \asymp  \P^w(\tau_m < T_{A[n/2,n]},
\tau_m < \bar{\eta}_n) = \P^w(\tau_m < \bar{\eta}_n)
\] 
We now show that (i) implies
\begin{equation}
\label{Eimport}
\P^z(\tau_m < \bar{\eta}_n)\leq \frac{c}{\log n}, \mbox{ for } z \in \CC_{4n}.
\end{equation}
Namely, by the same argument as in (i), applied to $n/2$ instead of $n$ and $m=n^4$ still, one gets
\[
\P^z(\tau_m<\bar{\eta}_n) \asymp \frac{1}{\log{n}} , \mbox{ for } z\in \CC_{2n}\setminus \CC_{3n/2}. 
\]
The uniform upper bound can easily be extended to all $z\in \CC_{2n}$ using strong Markov property and 
overshoot estimate (\ref{overshoot2}).
Now for $z \in \CC_{4n}\setminus \CC_{3n}$ we have
\[
\P^z(\tau_m<\bar{\eta}_n)=\P^z(\tau_m < \eta_n) + \P^z(\eta_n<\tau_m<\bar{\eta}_n),
\]
so that the upper bound in (i) together with strong Markov imply (\ref{Eimport}) for $z \in \CC_{4n}\setminus \CC_{3n}$. 
The remaining case $z \in \CC_{3n}\setminus \CC_{2n}$ is implied again by strong Markov inequality and an
overshoot estimate.
\ \ \ \endpf\\

Recall that we may assume $k\leq n/2$.
\begin{Lemma}
\label{LPesc}
If $ 0 \leq k \leq n/2$, and $j\in [k,n]_\kappa$ 
\[
\P^j(\tau_m <T_{[k,n]_\kappa}
) \leq \frac{c}{\sqrt{n}\log{n} }\left( \frac{1}{\psix{j-k+1}} +\frac{1}{\psix{n-j+1}}\right).
\]
\end{Lemma}
\Pf
Since $S_0=j$ the probability of not visiting $[k,n]_\kappa$ during interval $[1,\tau_{2n}]$ is bounded 
above by a constant times
\[
\frac{1}{\sqrt{n}\psix{j-k+1}} +\frac{1}{\sqrt{n}\psix{n-j+1}},
\]
due to Proposition \ref{prop2}(d),(e).

Now consider the first
time after $\tau_{2n}$ that the random walk either leaves $\CC_m$ or enters the disk
$\CC_n$.  
Estimate (\ref{Eimport}) says
that the probability of random walk leaving $\CC_m$ before entering
$\CC_n$ is bounded above by $c/\log n$.  
Hence one expects (also using an easy overshoot argument)
$O(\log n)$ ``excursions'' from $\CC_{4n}\setminus\CC_{2n}$ into $\CC_n$  before leaving $\CC_m$,
and for each such excursion there
is a positive probability, conditioned on the past of the walk, that the random walk visits $[k,n]_\kappa$
during that excursion due to Lemma \ref{gottahit}.
This gives the extra term $c/\log{n}$ in the above probability.\ \ \ \endpf

\section{Proof of Theorem \ref{Tmain}}
\label{PT}

Without loss of generality we assume $k,n \in \kappa \Z^+$
with $k \leq n/2$.   By Lemma \ref{loglem}(i), it suffices
to show that
\[  \P(\tau_m < T_{A[k,n]}) \leq \frac {c \sqrt k} {\sqrt n \, \log n} ,\]
where, as before, $m = n^4$.  
The above inequality will then imply 
\[
\P(\tau_{2n}<T_{A[k,n]})\leq \frac{\P(\tau_{m}<T_{A[k,n]})}{\P(\tau_m<T_{A[k,n]}|\tau_{2n}<T_{A[k,n]})}=\frac {c \sqrt k} {\sqrt n},
\]
since by Lemma \ref{loglem} (i)
\[
\P(\tau_m<T_{A[k,n]}|\tau_{2n}<T_{A[k,n]})\geq \frac{c}{\log{n}}.
\]

Let $T = T_{[k,n]_\kappa},
\hat T = T_{A[k,n]}, T^0 = T^0_{[k,n]_\kappa},
\hat T^0 = T_{A[k,n]}^0. $
Proposition \ref{prop2}(d),
Lemma \ref{loglem}(ii), and an easy overshoot estimate give
\[  \P(\tau_m < T) \leq  \frac {c \, \sqrt{k}} {\sqrt n \, \log n} . \]
Note that similarly we have
\[  \P(\tau_m < T) \geq  \frac {c \, \sqrt{k}} {\sqrt n \, \log n} , \]
with a different constant $c>0$, since Lemma \ref{loglem}(i) is two-sided bound, and 
the proof of Proposition \ref{prop2} parts (b) and (d) 
can be slightly modified to obtain two-sided bound of the same order as the upper bound.

We will show, in fact, that 
\[\P(\tau_m < \hat T)  - \P(\tau_m < T )
  \leq  \frac c{ \sqrt n \, \log n }  . \]
Note  that
$\P(\tau_m <\hat T)- \P(\tau_m < T)=
\P(T<\tau_m)-\P(\hat T <\tau_m )$
which equals, by (\ref{lastexit}) (note that under $\P$, $T^0=T$ and $\hat{T}^0=\hat{T}$),
\begin{eqnarray}
\lefteqn{\sum_{j\in [k,n]_\kappa} G_m(0,j)\P^j(\tau_m <T)-
\sum_{j\in [k,n]_\kappa} G_m(0,w_j)\P^{w_j}(\tau_m <\hat T)\;\;\;=}
    \hspace{.3in} \nonumber\\
& &\sum_{j\in [k,n]_\kappa} [G_m(0,j)- G_m(0,w_j)]\P^j(\tau_m <T)\;
+\label{Erea1}\\
& &\hspace{1in}
  \sum_{j\in [k,n]_\kappa} G_m(0,w_j)[
   \P^j(\tau_m <T)-\P^{w_j}(\tau_m <\hat T)].
\label{Erea2}
\end{eqnarray}
We will show that the sum in (\ref{Erea1}) is bounded above
 in absolute value
by $c/(\sqrt{n}\log{n})$ and  that the sum in (\ref{Erea2}) is 
bounded above by $c /(\sqrt{n}\log{n})$. We will not bound the
absolute value in (\ref{Erea2}).

Lemma \ref{Greenlemma} gives
\begin{equation}
\label{EGbd}
|G_m(0,j)-G_m(0,w_j)| \leq\frac{C}{j}.
\end{equation}
Lemma \ref{LPesc} gives
\[
\P^j(\tau_m <T) =O\left( \frac{1}{\sqrt{n} }
\left( \frac{1}{\psix{j-k+1}} +\frac{1}{\psix{n-j+1}}\right)\frac{1}{\log{n}}\right).
\]
The term in (\ref{Erea1}) is therefore bounded in absolute value by
\[
\sum_{j\in [k,n]_\kappa}\frac{C}{j}\left[ \frac{1}{\sqrt{n} }
\left( \frac{1}{\psix{j-k+1}} +\frac{1}{\psix{n-j+1}}\right)\frac{1}{\log{n}}\right]
  \leq \frac{c}{\sqrt n \; \log n}.\]
(Here and below we use the easy estimate:
\[  \sum_{j=-\infty}^\infty \frac{1}{|j-k| +1}\,  \frac{1}{(|j-l|+1)^{1/2}}
               \leq  2 \sum_{j=-\infty}^\infty \frac 1 {(|j| + 1)^{3/2}} < \infty.\;\;\;) \]

To estimate the term (\ref{Erea2}) define the function $f$ from $\L$ to $\R$ by
\[
f(z):=\sum_{j\in [k,n]_\kappa} G_m(z,w_j)
\; [\P^j(\tau_m < T) - \P^{w_j}(\tau_m <\hat T)],
\]
and note that (\ref{Erea2}) equals  $f({ 0})$.
Since $G_m(\cdot,w_j)$ is $p$-harmonic  on
$\{z:|z|<m\}\setminus \{w_j\}$, $f$ is $p$-harmonic on $\{z:|z|<m\}\setminus A[k,n]$,
and therefore it attains its maximum on  
$\{z:|z|\geq  m\}\cup A[k,n]$.
However, $f(z)=0$ for $z\geq m$, so  it suffices to show
\begin{equation}
\label{EsupA}
\max_{\ell \in [k,n]_\kappa} f(w_\ell) 
\leq \frac {c}{ \sqrt n \; \log n }.
\end{equation}
Fix $\ell \in [k,n]_\kappa $ and note by (\ref{lastexit})
\[
\sum_{j\in [k,n]_\kappa} G_m(w_\ell,w_j)\P^{w_j}(\tau_m <\hat T)=
    \P^{w_\ell}(\hat T^0<\tau_m)=1,
\]
and
\[
\sum_{j\in [k,n]_\kappa} G_m(\ell,j)\P^{j}(\tau_m <T)=
  \P^{\ell}(T^0<\tau_m)=1.\]
Hence,
\[ \sum_{j\in [k,n]_\kappa} 
G_m(w_\ell,w_j)\; [\P^j(\tau_m <T)-\P^{w_j}(\tau_m <\hat T)]=  \hspace{1in}\]
\begin{equation}
\label{Eagain} 
\hspace{1in} \sum_{j\in [k,n]_\kappa} [G_m(w_\ell,w_j)-G_m(\ell,j)]\;
 \P^j(\tau_m <T).
\end{equation}
Since $|w_\ell - w_j| \geq |\ell-j| - \kappa$,
Lemma \ref{Greenlemma} gives
\[
G_m(w_\ell,w_j)-G_m(\ell,j)\leq \frac{c}{|j-\ell| + 1}
\]
(note that we are not bounding the absolute value). 

Hence, 
\begin{eqnarray*}
\lefteqn{ \sum_{j\in[k,n]_\kappa } (G_m(w_\ell,w_j)-G_m(\ell,j)) \P^j(\tau_m <T)}
\nonumber\\
& \leq & \sum_{j\in [k,n]_\kappa}  \frac{1}{\sqrt{n}}
\;  \frac{1}{|j-\ell| + 1} \, 
 \left( \frac{1}{\psix{j-k+1}} +\frac{1}{\psix{n-j+1}}\right)\frac{1}{\log{n}} \\
& \leq & \frac{c} {\sqrt n  \log n}. \ \ \mbox{ \endpf}
\end{eqnarray*}

\noindent
{\bf Acknowledgment.} We are grateful to Yasunari Fukai for useful conversations.

\end{document}